\documentstyle{article}[15pt]
\input amssym.def

\catcode`@=12
 \long\def\@makefntext#1{\noindent #1}
\newskip\tabcentering \tabcentering=1000pt plus 1000pt minus 1000pt

\def\MCH#1#2{\setbox0=\hbox{\raise#1\hbox{#2}}\smash{\box0}}

\def\@evenfoot{}\def\@oddfoot{}

\def\sec#1{\vspace{5mm}\leftline{\bf #1}\vspace{3mm}}

\catcode`@=12

\def\bc{\begin{center}}
\def\ec{\end{center}}

\def\hang{\hangindent\parindent}
\def\textindent#1{\indent\llap{\qquad #1\ \ \enspace}\ignorespaces}
\def\ref{\par\hang\textindent}

\def\a1{(a_1, a_2, \cdots, a_n)}

\def\a{\alpha}

\begin{document}
\thispagestyle{empty}
\vspace*{-3.0truecm}
\noindent
\vspace{1 true cm}
 \bc{\large\bf Delta  shocks   and  vacuum   states for  the isentropic  magnetogasdynamics   equations
 for  Chaplygin gas  as pressure and magnetic field vanish $^{ **}$

\footnotetext{$^{*}$Corresponding author. Tel: +86-0591-83852790.\\
\indent \,\,\,\,\,\,\,\,E-mail address:zqshao@fzu.edu.cn.\\
\indent \,\,\,\,\,$^{**}$Supported by the National Natural Science
Foundation of China (No. 70371025),  the Scientific Research
Foundation of the
 Ministry of Education of China (No. 02JA790014),
  the Natural Science Foundation of Fujian Province of China   (No.
 2012J01006)   and  the Science and Technology Developmental Foundation
of Fuzhou University (No. 2004-XQ-16).  }}\ec

 \vspace*{0.2 true cm}
\bc{\bf  Zhiqiang  Shao$^{a, *}$    \\
{\it $^{a}$Department of Mathematics,  Fuzhou University,  Fuzhou 350002,
China}
 }\ec

 \vspace*{2.5 true mm}
\setlength{\unitlength}{1cm}
\begin{picture}(20,0.1)
\put(-0.6,0){\line(1,0){14.5}}
\end{picture}

 \vspace*{2.5 true mm}
\noindent{\small {\small\bf Abstract}

 \vspace*{2.5 true mm}This paper is concerned with the Riemann problem for  the isentropic Chaplygin gas magnetogasdynamics   equations
   and the formation of delta shocks and vacuum states as pressure and magnetic field vanish.
Firstly, the Riemann problem of   the isentropic  magnetogasdynamics   equations
 for  Chaplygin gas is solved analytically. Secondly, it is rigorously proved that, as both the pressure and the magnetic field vanish,   the Riemann solution containing two shock waves tends to
 a  delta
 shock solution to  the transport equations,
 and the intermediate density between the two shocks tends to a weighted $\delta$-measure which forms the  delta
 shock; while the Riemann solution   containing two rarefaction waves 
 tends to a two-contact-discontinuity solution to the transport equations, the  termediate state between the two contact discontinuities
is a vacuum state.

 \vspace*{2.5 true mm}
\noindent{\small {\small\bf MSC: } 35L65;  35L67

 \vspace*{2.5 true mm}
\noindent{\small {\small\bf Keywords:} Isentropic  magnetogasdynamics;   Chaplygin gas;   Riemann problem; Transport equations;  Delta shock wave;   Vacuum state

 \vspace*{2.5 true mm}
\setlength{\unitlength}{1cm}
\begin{picture}(20,0.1)
\put(-0.6,0){\line(1,0){14.5}}
\end{picture}



\baselineskip 15pt
 \sec{\Large\bf 1.\quad  Introduction }In this paper, we are concerned with the   system  of  conservation  law   governing  the one-dimensional unsteady simple flow of an isentropic, inviscid and perfectly conducting compressible fluid subjected to a transverse magnetic field    (see [10,  11]):
$$ \left\{\begin{array}{ll} \rho_{t}+(\rho u)_{x}=0,\\(\rho u)_{t}+
(p+\rho u^{2}+B^{2}/2\mu)  _{x}=0,\end{array}\right .\eqno{(1.1)}
$$
where $\rho>0$,      $u,$  $p$, $B$   and  $\mu>0$  represent the density, velocity, pressure, transverse magnetic field and magnetic
permeability, respectively; $p$ and $B$ are known functions defined as $$p=-\frac{k_{1}}{\rho}  \,\,\,\,\,\,\,\,\,\eqno{(1.2)}
$$  and $B=k_{2}\rho,$ where $k_{1}$ and $k_{2}$ are positive
constants. The independent variables $t$ and
$x$ denote time and space,  respectively.  The adiabatic exponent in (1.2) can be viewed as $\gamma=-1$ by comparing
with the state equation $ p=k_{1}\rho^{\gamma}   $  with $\gamma\geq 1$   for the polytropic gas.
The gas (1.2) whose adiabatic constant $\gamma=-1$
 is  usually called as the  Chaplygin gas.

 For the  isentropic Chaplygin gas  Euler equations, Brenier [1] firstly studied the 1-D Riemann problem and obtained solutions
with concentration when initial data belong to a certain domain in the phase plane. Furthermore, Guo,
Sheng and Zhang [6] abandoned this constrain and constructively obtained the global solutions to the
1-D Riemann problem, in which the $\delta$-shock developed. Moreover, they also systematically studied the
2-D Riemann problem for isentropic Chaplygin gas equations. For the 2-D case, we can also refer to
[9] in which D. Serre studied the interaction of the pressure waves for the 2-D isentropic irrotational
Chaplygin gas and constructively proved the existence of transonic solutions for two cases, saddle and
vortex of 2-D Riemann problem. Recently,
Sheng, Wang and Yin [13] and Wang [15] studied the Riemann problem  for the generalized Chaplygin gas  and obtained the solutions to the
Riemann problem and the interactions of elementary waves. The Riemann solutions to the transport equations in zero-pressure
flow in gas dynamics were presented by Sheng and Zhang in [14], in which delta shocks and vacuum states appeared.

In related researchs of the $\delta$-shock  waves,  one very important and interesting  topic is to study the
phenomena of concentration and cavitation and the formation of
$\delta$-shock  waves   and  vacuum   states in solutions.  In earlier paper [4],  Chen and Liu [4] studied  the formation of
$\delta$-shocks   and  vacuum   states of the  Riemann solutions to the  isentropic  Euler equations  for polytropic gas  as $\varepsilon\rightarrow 0$,   in which they took  the equation of state as $ P=\varepsilon p  $  for $ p=\rho^{\gamma}/\gamma   $ $(\gamma> 1).$  Further,  they also obtained the same results for the Euler equations   for nonisentropic fluids in [5].
 The same problem for the the  isentropic Euler equations  for isothermal case was studied by Li [7],  in which he proved that  when  temperature drops to zero,  the solution containing two shock waves converges to the  delta
 shock solution to  the transport equations  and the solution   containing two rarefaction waves converges to the
 solution involving vacuum to the transport equations.   Then, the results were extended to the relativistic Euler
 equations for polytropic gas by Yin and Sheng [17]  and for Chaplygin gas by Yin and Song [18], the  isentropic  Euler equations  for  the generalized Chaplygin gas by
Sheng, Wang and Yin [13] and for modified Chaplygin gas  by Yang  and  Wang [16], the perturbed Aw-Rascle
model by Shen  and Sun [12],     the isentropic magnetogasdynamics equations for polytropic gas by  Shen [11],
  the   generalized pressureless gas dynamics  model with a scaled pressure term by Mitrovic
and Nedeljkov [8],  etc.

In this paper, we study the Riemann problem of  the isentropic  magnetogasdynamics   equations
 for  Chaplygin gas 
   and the formation of delta shocks and vacuum states as pressure and magnetic field vanish. The organization of the paper is as follows. In Sections 2 and 3,  the Riemann problems for  the isentropic  Chaplygin gas   magnetogasdynamics   equations
      and the transport equations are analyzed by characteristic analysis. In Sections 4 and 5, we investigate the formation of
$\delta$-shocks   and  vacuum   states of the  Riemann solutions to  the isentropic  magnetogasdynamics   equations
 for  Chaplygin gas  as pressure and magnetic field vanish.

\baselineskip 15pt
 \sec{\Large\bf 2.\quad    Riemann problem for
  system
 (1.1)-(1.2) }In this section, we discuss the Riemann solutions of (1.1) and (1.2) with initial data$$ (\rho, u)( x, 0) =(\rho_{\pm},
u_{\pm}),\,\,\,\,\,\,\,\,\,\,\,\,\pm x> 0,\eqno{(2.1)}$$
where $\rho_{\pm} >0$  and $u_{\pm}$  are arbitrary constants.

For smooth solution, system (1.1) is equivalent to
$$\left(\begin{array}{cccc}\rho\\u
 \end{array}\right)_{t}+
\left(\begin{array}{cccc}u &\rho\\w^{2}/\rho &u
 \end{array}\right)\left(\begin{array}{cccc}\rho\\u
 \end{array}\right)_{x}=0,\eqno{(2.2)}$$
where $w=(c^{2}+b^{2})^{1/2}$ is the magneto-acoustic speed with $c=(p'(\rho))^{1/2} $ as the local sound speed and $b=(B^{2}(\rho)/\mu\rho)^{1/2}$ the Alfven speed. Here, prime denotes differentiation with respect to $\rho$.
 The eigenvalues of system (1.1) and (1.2) are
$$ \lambda_{1}=u-\sqrt{\frac{k_{1}}{\rho^{2}}+\frac{k_{2}^{2}\rho}{\mu}},\,\,\,\,\,\,
\,\lambda_{2}=u+\sqrt{\frac{k_{1}}{\rho^{2}}+\frac{k_{2}^{2}\rho}{\mu}}.\eqno{}
$$
Therefore,  system (1.1) and (1.2) is strictly
hyperbolic for $\rho> 0$.

The corresponding right eigenvectors   are
$$ \overrightarrow{r_{1}}=(-\rho, \sqrt{\frac{k_{1}}{\rho^{2}}+\frac{k_{2}^{2}\rho}{\mu}})^{T},\,\,\,\,\,\,\overrightarrow{r_{2}}=(\rho, \sqrt{\frac{k_{1}}{\rho^{2}}+\frac{k_{2}^{2}\rho}{\mu}})^{T}.\eqno{}
$$
By simple calculation, we get
$$ \nabla\lambda_{i}\cdot \overrightarrow{r_{i}}=\frac{3k_{2}^{2}\rho}{2\mu\sqrt{\frac{k_{1}}{\rho^{2}}+\frac{k_{2}^{2}\rho}{\mu}}}\neq 0,\,\,\,\,\,i=1, 2.\eqno{}
   $$
Therefore, both the characteristic fields
 are genuinely nonlinear.

Since system (1.1), (1.2) and the Riemann data (2.1) are invariant under stretching of coordinates:
$(x,t)\rightarrow (\alpha x,\alpha t)~(\alpha$ is constant),  we seek the self-similar solution $$(\rho,u)(x, t)=(\rho,u)(\xi),\,\,\,\,\xi=\frac{x}{t}.$$
 Then  Riemann problem (1.1), (1.2) and (2.1) is reduced to the following boundary value problem of  ordinary differential equations:
$$ \left\{\begin{array}{ll}
   -\xi\rho_{\xi}+(\rho u)_{\xi}=0,\\
   -\xi(\rho u)_{\xi}+\big(-\frac{k_{1}}{\rho}
   +\rho u^{2}+\frac{(k_{2}\rho)^{2}}{2\mu}\big)_{\xi}=0,\end{array}\right .\eqno{(2.3)}
$$
with $(\rho,u)(\pm\infty)=(\rho_{\pm},u_{\pm}).$\\
\indent
For any smooth solution, system (2.3) can be written as
$$
\left(\begin{array}{cccc}u-\xi &\rho\\-\xi u+\frac{k_{1}}{\rho^{2}}+u^{2}+\frac{k_{2}^{2}\rho}{\mu} &-\xi \rho+2\rho u
 \end{array}\right)\left(\begin{array}{cccc}\rho_{\xi}\\u_{\xi}
 \end{array}\right)=0.\eqno{(2.4)}$$
It provides either general solutions (constant states)
$$(\rho,u)(\xi)={\mathrm constant }   \,\,\,\,\,(\rho> 0)$$
or  singular solutions called   the rarefaction waves $R_{1}$  and $R_{2}$
which denote, respectively, 1-rarefaction waves and 2-rarefaction waves,
 $$R_{1}:
\left\{
  \begin{array}{ll}
    \xi=\lambda_{1}=u-\sqrt{\frac{k_{1}}{\rho^{2}}+\frac{k_{2}^{2}\rho}{\mu}}, \\
    u-u_{-}=-\int^{\rho}_{\rho_{-}}\frac{\sqrt{\frac{k_{1}}{s^{2}}+\frac{k_{2}^{2}s}{\mu}}}{s}ds,\,\,\,\,
    \rho<\rho_{-},
  \end{array}
\right.
\eqno{(2.5)}  $$
and$$R_{2}:
\left\{
  \begin{array}{ll}
    \xi=\lambda_{2}=u+\sqrt{\frac{k_{1}}{\rho^{2}}+\frac{k_{2}^{2}\rho}{\mu}}, \\
    u-u_{-}=\int^{\rho}_{\rho_{-}}\frac{\sqrt{\frac{k_{1}}{s^{2}}+\frac{k_{2}^{2}s}{\mu}}}{s}ds,\,\,\,\,
    \rho>\rho_{-}.
  \end{array}
\right.
\eqno{(2.6)}  $$

 Differentiating the second equation of (2.5) with respect to $\rho$ yields $u_{\rho}=-\frac{\sqrt{\frac{k_{1}}{\rho^{2}}+\frac{k_{2}^{2}\rho}{\mu}}}{\rho}<0,$  and subsequently, $$u_{\rho\rho}=\frac{\frac{4k_{1}}{\rho^{2}}+\frac{k_{2}^{2}\rho}{\mu}}{2\rho^{2}\sqrt{\frac{k_{1}}{\rho^{2}}
+\frac{k_{2}^{2}\rho}{\mu}}}>0,$$    which mean that the 1-rarefaction wave  curve $R_{1}$ is  monotonic decreasing  and convex  in the $(\rho, u)$ plane $(\rho>0 )$.  Similarly, from  the second equation of (2.6), we have $u_{\rho}>0$   and   $u_{\rho\rho}<0$,   which mean that the 2-rarefaction wave  curve $R_{2}$ is    monotonic increasing  and   concave  in the $(\rho, u)$ plane $(\rho>0 )$. In addition,   it can be  verified that $\lim\limits_{\rho\rightarrow 0^{+}}u=+\infty$  for  the 1-rarefaction wave  curve $R_{1}$, which implies that
$R_{1}$ has the $u$-axis as its  asymptotic line.
It can also be proved that $\lim\limits_{\rho\rightarrow  +\infty}u=+\infty$  for  the  2-rarefaction wave  curve $R_{2}$.

\indent
For a bounded discontinuity at $\xi=\sigma,$ the Rankine-Hugoniot conditions hold:
$$\left\{
    \begin{array}{ll}
      -\sigma[\rho]+[\rho u]=0, \\
     -\sigma[\rho u]+[-\frac{k_{1}}{\rho}
   +\rho u^{2}+\frac{(k_{2}\rho)^{2}}{2\mu}]=0,
    \end{array}
  \right.\eqno (2.7)
$$
where $[\rho]=\rho -\rho_{-},$   etc.  Solving (2.7), we obtain two   shock waves $S_{1}$ and $S_{2}$
$$S_{1}: \left\{
             \begin{array}{ll}
               \sigma=u_{-}-\rho\sqrt{\frac{1}{\rho\rho_{-}}\Big(\frac{k_{1}}{\rho\rho_{-}}+\frac{k_{2}^{2}(\rho+\rho_{-})}{2\mu}\Big)}, \\
             u=u_{-}-\sqrt{\frac{1}{\rho\rho_{-}}\Big(\frac{k_{1}}{\rho\rho_{-}}
             +\frac{k_{2}^{2}(\rho+\rho_{-})}{2\mu}\Big)}(\rho-\rho_{-}),  \,\,\,
              \rho>\rho_{-},
             \end{array}
           \right.
\eqno{(2.8)}  $$
$$S_{2}:
\left\{
             \begin{array}{ll}
               \sigma=u_{-}+\rho\sqrt{\frac{1}{\rho\rho_{-}}\Big(\frac{k_{1}}{\rho\rho_{-}}+\frac{k_{2}^{2}(\rho+\rho_{-})}{2\mu}\Big)}, \\
             u=u_{-}+\sqrt{\frac{1}{\rho\rho_{-}}\Big(\frac{k_{1}}{\rho\rho_{-}}
             +\frac{k_{2}^{2}(\rho+\rho_{-})}{2\mu}\Big)}(\rho-\rho_{-}),  \,\,\,
    \rho<\rho_{-}.
  \end{array}
\right.
\eqno{(2.9)}  $$

Differentiating the second equation of (2.8) with respect to $\rho$,   for  $\rho>\rho_{-}$ we have $$u_{\rho}=-\frac{1}{2\sqrt{\frac{1}{\rho\rho_{-}}\Big(\frac{k_{1}}{\rho\rho_{-}}
             +\frac{k_{2}^{2}(\rho+\rho_{-})}{2\mu}\Big)}}\Big(\frac{2k_{1}}{\rho_{-}\rho^{3}}
             +\frac{k_{2}^{2}}{\rho_{-}\mu}+\frac{k_{2}^{2}}{2\rho\mu}+\frac{k_{2}^{2}\rho_{-}}{2\rho^{2}\mu}\Big)<0,   $$     which mean that the 1-shock  curve $S_{1}$ is  monotonic decreasing   in the $(\rho, u)$ plane $(\rho>0 )$.  Similarly, from  the second equation of (2.9), for  $\rho<\rho_{-}$ we have $u_{\rho}>0,$   which mean that the 2-shock wave  curve $S_{2}$ is    monotonic increasing    in the $(\rho, u)$ plane $(\rho>0 )$.   In addition,   it is easily derived  from (2.9) that  $\lim\limits_{\rho\rightarrow  0^{+}}u=-\infty$  for  the 2-shock  curve $S_{2},$ which implies that
$S_{2}$ has the $u$-axis as its  asymptotic line. It can also be derived  from (2.8)  that $\lim\limits_{\rho\rightarrow+\infty}u=-\infty$  for  the 1-shock  curve $S_{1}$.

In the phase plane $(\rho>0,\, u \in\mathbf{R})$, through  point $(\rho_{-},u_{-})$, we draw the elementary  wave curves  $R_{1}$, $R_{2}$, $S_{1}$ and $S_{2}$,  respectively. Then   the phase plane is divided into four   regions I, II, III  and  IV$(\rho_{-}, u_{-})$    (see Fig. 1).

By  the analysis  method in phase plane, for any  given state $(\rho_{+}, u_{+}),$
 one can construct the  Riemann solutions   as follows:

 (1) $(\rho_{+}, u_{+})\in I(\rho_{-}, u_{-}): R_{1}+R_{2};$

 (2) $(\rho_{+}, u_{+})\in II(\rho_{-}, u_{-}): R_{1}+S_{2};$

 (3) $(\rho_{+}, u_{+})\in III(\rho_{-}, u_{-}): S_{1}+R_{2};$

 (4) $(\rho_{+}, u_{+})\in IV(\rho_{-}, u_{-}): S_{1}+S_{2}.$

\hspace{65mm}\setlength{\unitlength}{0.8mm}\begin{picture}(80,66)
\put(-50,2){\vector(0,2){50}}
 \put(-48,0){\vector(2,0){103}}  \put(-53,49){$\rho$}
\put(56,-1){$u$}
\put(-38,7){$S_{2}$}\put(6,5){}
\put(-32,48){$S_{1}$}\put(32,48){$R_{2}$}\put(42,7){$R_{1}$}
\put(5,19){$(\rho_{-},
u_{-})$}
\put(-3,6){II }\put(44,6){}\put(-3,29){III}
\put(-36,20){IV}\put(36,20){I}
\qbezier(50,4)(-14,12)(-42,52)\qbezier(-50,4)(14,12)(42,52)
\end{picture}
\vspace{0.6mm}  \vskip 0.2in \centerline{\bf Fig. 1.\,\,    Curves of  elementary waves.
   } \vskip 0.1in \indent

Thus  we  have proved the following result
\vskip 0.1in
\noindent{\small {\small\bf Theorem 1.}   For
Riemann problem (1.1),  (1.2)  and (2.1),   there exists   a unique entropy solution,  which  consists of    shock waves, rarefaction waves, and   constant states.
 $$$$
 $$$$

\baselineskip 15pt
 \sec{\Large\bf 3.\quad   Riemann problem for the transport equations }The Riemann solutions to the transport equations in zero-pressure flow were presented by Sheng and Zhang in [14]. The
Riemann problem to the transport equations are
$$ \left\{\begin{array}{ll} \rho_{t}+(\rho u)_{x}=0,\\(\rho u)_{t}+
 (\rho u^{2} )_{x}=0\end{array}\right .\eqno{(3.1)}
$$
  with initial data$$ (\rho,u)( x, 0) =(\rho_{\pm},
u_{\pm}),\,\,\,\,\,\,\,\,\,\,\,\,\pm x> 0.\eqno{(3.2)}$$

The system
 has  a double
eigenvalue
$$ \lambda=u\eqno{}
$$
and only one  right eigenvector
$$ \overrightarrow{r}=(r, 0)^{T}.\eqno{}
$$
By a direct calculation,
$$ \nabla\lambda\cdot \overrightarrow{r}\equiv 0.\eqno{}
   $$
Thus (3.1) is nonstrictly hyperbolic and  $\lambda$  is
 linearly degenerate.

As usual,   we seek the self-similar solution $$(\rho,u)(x, t)=(\rho,u)(\xi),\,\,\,\,\xi=\frac{x}{t}.$$
 Then  Riemann problem (3.1) and (3.2) is reduced to the following boundary value problem of  ordinary differential equations:
$$ \left\{\begin{array}{ll}
   -\xi \rho_{\xi}+(\rho u)_{\xi}=0,\\
   -\xi (\rho u)_{\xi}+(\rho u^{2})_{\xi}=0,\end{array}\right .\eqno{(3.3)}
$$
with $(\rho,u)(\pm\infty)=(\rho_{\pm},u_{\pm}).$\\
\indent
For any smooth solution, system (3.3) can be written as
$$
\left(\begin{array}{cccc}u-\xi &\rho\\0 &\rho(u-\xi)
 \end{array}\right)\left(\begin{array}{cccc}\rho_{\xi}\\u_{\xi}
 \end{array}\right)=0.\eqno{(3.4)}$$
It provides either the general solution (constant state)
$$(\rho,u)(\xi)={\mathrm constant }   \,\,\,\,\,(\rho\neq 0)$$
or   the singular solution
$$
\left\{
  \begin{array}{ll}
    \rho=0,\\u=\xi,
  \end{array}
\right.
\eqno{(3.5)}  $$
which is called   the vacuum state (see [14]), where $u(\xi)$ is an arbitrary smooth function.

\indent
For a bounded discontinuity at $\xi=\sigma,$ the Rankine-Hugoniot condition holds:
$$\left\{
    \begin{array}{ll}
      -\sigma[\rho]+[\rho u]=0, \\
     -\sigma[\rho u]+[\rho u^{2}]=0,
    \end{array}
  \right.\eqno (3.6)
$$
 where $[q]=q_{+} -q_{-}$ denotes the jump of $q$  across the discontinuity. By solving (3.6), we obtain   $$J : \xi = \sigma= u_{-}(=\lambda_{-}) =u_{+}(=\lambda_{+}), \eqno (3.7)
$$
which is a contact discontinuity. It is a slip line and  just the characteristic of solutions on both its sides in $(x, t)$-plane.

The Riemann problem (3.1) and (3.2) can be solved by contact discontinuities, vacuum  or delta shock wave connecting two constant states $(u_{\pm},v_{\pm})$.

For the case $u_{-}< u_{+}$, there is no characteristic passing through the region $u_{-}t < x < u_{+}t $  and  the vacuum  appears in this  region. The solution can be expressed as
$$ (\rho,u)(\xi)=\left\{\begin{array}{ll} (\rho_{-}, u_{-}),\,\,\,\,\,\,-\infty<x<u_{-}, \\(0,\xi),\,\,\,\,\,\,\,\,\,\,\,\,\,\,\,\,\,u_{-} \leq \xi \leq u_{+},\\(\rho_{+},u_{+}), \,\,\,\,\,\,u_{+}<\xi < +\infty.
\end{array}\right .\eqno{(3.8)}
$$

For the case $u_{-}= u_{+}$, it is easy to see that the constant states $(\rho_{\pm},u_{\pm})$ can be connected by a contact discontinuity.

\hspace{65mm}\setlength{\unitlength}{0.8mm}\begin{picture}(80,66)
\put(0,0){\line(-4,5){25}}
 \put(0,0){\line(1,3){12}}
\put(8,42){$\frac{x}{t}=u_{+}$}
\put(15,10){$(\rho_{+},
u_{+})$}
\put(-36,8){$(\rho_{-},
u_{-})$}

\put(-3,16){$\Omega$ }
\put(-34,34){$\frac{x}{t}=u_{-}$}

\qbezier(3,10)(-2,13)(-5,7)

 \put(0,0){\vector(0,2){48}}
\put(-45,0){\vector(2,0){98}} \put(-3,-4){$O$} \put(-4,45){$t$}
\put(54,-1){$x$}
\end{picture}
\vspace{0.6mm}  \vskip 0.2in \centerline{\bf Fig. 2.\, Characteristics overlap domain.} \vskip 0.1in \indent

For the case $u_{-}> u_{+}$,   the characteristic lines  originating from the origin  will overlap
in a domain $\Omega$,  as shown in Fig. 2. So, singularity must happen in $\Omega$. It is easy to know that the singularity is impossible to
be a jump with finite amplitude because the Rankine-Hugoniot condition is not satisfied on the bounded jump. In other
words, there is no solution which is piecewise smooth and bounded. Motivated by  [14], we seek solutions with
delta distribution at the jump.

To do so, a
  two-dimensional weighted delta function $w(s)\delta_{L}$ supported on a smooth curve $ L=\{(t(s),x(s)):a<s<b\}$ is defined by
$$ \langle w(s)\delta_{L},\varphi\rangle=\int_{a}^{b}w(s)\varphi(t(s),x(s))ds \eqno{(3.9)}$$
for any  $\varphi\in C^{\infty}_{0}(R\times R_{+}).$ \\

\indent
Let us consider a  solution of  (3.1) and (3.2) of the form
$$(\rho,u)(x, t)=\left\{
                  \begin{array}{ll}
                    (\rho_{-}, u_{-}), & \hbox{$x<\sigma t $,} \\
                    (w(t)\delta(x-\sigma t),\sigma), & \hbox{$x=\sigma t$,} \\
                    (\rho_{+}, u_{+}), & \hbox{$x>\sigma t$,}
                  \end{array}
                \right.
\eqno{(3.10)}$$
where  $\sigma$  is a   constant,   $w(t)\in C^{1}[0, +\infty)$, and $\delta(\cdot)$ is the standard Dirac measure. $x(t)$, $w(t)$ and $\sigma$  are  the location,  weight and  velocity  of the delta shock,
 respectively. Then the following generalized Rankine-Hugoniot conditions hold:
$$
\left\{
     \begin{array}{ll}
       \frac{dx(t)}{dt}=\sigma, \\
       \frac{dw(t)}{dt}=\sigma [\rho]-[\rho u], \\
       \frac{d(w(t)\sigma)}{dt}=\sigma [\rho u]-[\rho u^{2}],
     \end{array}
   \right.\eqno (3.11)
$$
where $[\rho]= \rho_{+}-\rho_{-}$, with initial data
$$(x, w)(0) = (0,  0).\eqno{(3.12)}$$

\indent
In addition, to guarantee uniqueness,  the delta  shock wave should satisfy the entropy
condition:
$$u_{+}<\sigma<u_{-}.$$

 Solving the system of  simple ordinary differential equations (3.11) with initial data
 (3.12),
we have, when $\rho_{-}=\rho_{+}$,$$x(t)=\frac{1}{2}(u_{-}+u_{+})t,~~w(t)=(\rho_{-}u_{-} -\rho_{+}u_{+})t,$$$$\sigma=\frac{1}{2}(u_{-}+u_{+});$$
 when $\rho_{-}\neq  \rho_{+}$,$$ x(t)=\frac{\sqrt{\rho_{-}}u_{-}+\sqrt{\rho_{+}}u_{+}}{\sqrt{\rho_{-}}+\sqrt{\rho_{+}}}t,~~w(t)=\sqrt{\rho_{-}\rho_{+}}(u_{-}-u_{+})t,$$
 $$\sigma=\frac{\sqrt{\rho_{-}}u_{-}+\sqrt{\rho_{+}}u_{+}}{\sqrt{\rho_{-}}+\sqrt{\rho_{+}}}.$$

\baselineskip 15pt
 \sec{\Large\bf 4.\quad  Formation of $\delta$-shocks}
In this section,  we study the formation of $\delta$-shock waves in the  Riemann solutions of system  (1.1) and (1.2) in the case  $(\rho_{+}, u_{+})\in IV(\rho_{-}, u_{-})$  with  $u_{-}> u_{+}$  as  both the pressure and the magnetic field vanish.

\baselineskip 15pt
\sec{\large\bf 4.1.\quad    Limit behavior of  Riemann solutions as $k_{1}, k_{2}\rightarrow 0$
}
When $(\rho_{+}, u_{+})\in IV(\rho_{-}, u_{-}), $  for each pair of fixed $k_{1}> 0$  and
$k_{2}>0$, suppose that  $(\rho_{\ast}, u_{\ast}) $ is the intermediate state   connected    with $(\rho_{-}, u_{-})$ by a 1-shock $S_{1}$
with speed $\sigma_{1}$  and $(\rho_{+}, u_{+})$   by a 2-shock $S_{2}$
with speed $\sigma_{2}$. Then it follows
$$S_{1}: \left\{
             \begin{array}{ll}
               \sigma_{1}=u_{-}-\rho_{\ast}\sqrt{\frac{1}{\rho_{\ast}\rho_{-}}
               \Big(\frac{k_{1}}{\rho_{\ast}\rho_{-}}+\frac{k_{2}^{2}(\rho_{\ast}+\rho_{-})}{2\mu}\Big)}, \\
             u_{\ast}=u_{-}-\sqrt{\frac{1}{\rho_{\ast}\rho_{-}}\Big(\frac{k_{1}}{\rho_{\ast}\rho_{-}}
             +\frac{k_{2}^{2}(\rho_{\ast}+\rho_{-})}{2\mu}\Big)}(\rho_{\ast}-\rho_{-}),  \,\,\,
              \rho_{\ast}>\rho_{-},
             \end{array}
           \right.
\eqno{(4.1)}  $$
$$S_{2}:
\left\{
             \begin{array}{ll}
               \sigma_{2}=u_{\ast}+\rho_{+}\sqrt{\frac{1}{\rho_{+}\rho_{\ast}}\Big(\frac{k_{1}}
               {\rho_{+}\rho_{\ast}}+\frac{k_{2}^{2}(\rho_{+}+\rho_{\ast})}{2\mu}\Big)}, \\
             u_{+}=u_{\ast}+\sqrt{\frac{1}{\rho_{+}\rho_{\ast}}\Big(\frac{k_{1}}{\rho_{+}\rho_{\ast}}
             +\frac{k_{2}^{2}(\rho_{+}+\rho_{\ast})}{2\mu}\Big)}(\rho_{+}-\rho_{\ast}),  \,\,\,
    \rho_{+}<\rho_{\ast}.
  \end{array}
\right.
\eqno{(4.2)}  $$
The addition of (4.1) and (4.2) gives
$$\hspace{-35mm}u_{-}-u_{+}
             =\sqrt{\frac{1}{\rho_{-}}-\frac{1}{\rho_{\ast}}}
             \sqrt{k_{1}(\frac{1}{\rho_{-}}-\frac{1}{\rho_{\ast}})
             +\frac{k_{2}^{2}(\rho_{\ast}^{2}-\rho_{-}^{2})}{2\mu}}$$$$
             +\sqrt{\frac{1}{\rho_{+}}-\frac{1}{\rho_{\ast}}}
             \sqrt{k_{1}(\frac{1}{\rho_{+}}-\frac{1}{\rho_{\ast}})
             +\frac{k_{2}^{2}(\rho_{\ast}^{2}-\rho_{+}^{2})}{2\mu}}
             ,\,\,\,\,\,\rho_{\ast}>\rho_{\pm}.\eqno{(4.3)}$$
For any given $\rho_{\pm}>0$, if $\lim\limits_{k_{1}, k_{2}\rightarrow 0}\rho_{\ast}=M\in[\max(\rho_{-}, \rho_{+}), +\infty),$  then by taking the limit $k_{1}, k_{2}\rightarrow 0$  in (4.3), we have $u_{-}-u_{+}=0,$
which contradicts with $u_{-}>u_{+}.$
Therefore, $\lim\limits_{k_{1}, k_{2}\rightarrow 0}\rho_{\ast}=+\infty.
$  Letting $k_{1}, k_{2}\rightarrow 0$ in (4.3),  we obtain  the following result.

\noindent{\small \small\bf Lemma 1.}
$$\lim\limits_{k_{1}, k_{2}\rightarrow 0}k_{2}^{2}\rho_{\ast}^{2}=\frac{2\mu\rho_{-}\rho_{+}(u_{-}-u_{+})^{2}}{(\sqrt{\rho_{-}}+\sqrt{\rho_{+}})^{2}}.
\eqno{(4.4)}
$$

\noindent{\small \small\bf Lemma 2.}
$$\lim\limits_{k_{1}, k_{2}\rightarrow 0}u_{\ast}=\lim\limits_{k_{1}, k_{2}\rightarrow 0}\sigma_{1}=\lim\limits_{k_{1}, k_{2}\rightarrow 0}\sigma_{2}=\frac{\sqrt{\rho_{-}}u_{-}+\sqrt{\rho_{+}}u_{+}}{\sqrt{\rho_{-}}+\sqrt{\rho_{+}}}=\sigma,
\eqno{(4.5)}
$$
$$\lim\limits_{k_{1}, k_{2}\rightarrow 0}\int^{\sigma_{2}t}_{\sigma_{1}t}\rho_{\ast}dx=(\sigma[\rho]-[\rho u])t
=\sqrt{\rho_{-}\rho_{+}}(u_{-}-u_{+})t=w(t).
\eqno{(4.6)}
$$
\vskip 0.1in
\noindent{\small\small\bf Proof.}  Letting $k_{1}, k_{2}\rightarrow 0$ in (4.1)  and noting Lemma 4.1, we have
$$\hspace{-18mm}\lim\limits_{k_{1}, k_{2}\rightarrow 0}u_{\ast}=u_{-}-\lim\limits_{k_{1}, k_{2}\rightarrow 0} \sqrt{\frac{1}{\rho_{-}}-\frac{1}{\rho_{\ast}}}
             \sqrt{k_{1}(\frac{1}{\rho_{-}}-\frac{1}{\rho_{\ast}})
             +\frac{k_{2}^{2}(\rho_{\ast}^{2}-\rho_{-}^{2})}{2\mu}}$$
             $$=u_{-}-\sqrt{\frac{1}{\rho_{-}}}\sqrt{\frac{\rho_{-}\rho_{+}(u_{-}-u_{+})^{2}}{(\sqrt{\rho_{-}}+\sqrt{\rho_{+}})^{2}}}
             =\frac{\sqrt{\rho_{-}}u_{-}+\sqrt{\rho_{+}}u_{+}}{\sqrt{\rho_{-}}+\sqrt{\rho_{+}}}=\sigma.
\eqno{(4.7)}
$$
From  the first equation of (4.1),  by Lemma 4.1,   we obtain
$$\hspace{-18mm}\lim\limits_{k_{1}, k_{2}\rightarrow 0}\sigma_{1}=u_{-}-\lim\limits_{k_{1}, k_{2}\rightarrow 0}
             \sqrt{\frac{k_{1}}{\rho^{2}_{-}}
             +\frac{k_{2}^{2}\rho_{\ast}^{2}(\frac{1}{\rho_{-}}+\frac{1}{\rho_{\ast}})}{2\mu}}$$
             $$=u_{-}-\sqrt{\frac{\rho_{+}(u_{-}-u_{+})^{2}}{(\sqrt{\rho_{-}}+\sqrt{\rho_{+}})^{2}}}
             =\frac{\sqrt{\rho_{-}}u_{-}+\sqrt{\rho_{+}}u_{+}}{\sqrt{\rho_{-}}+\sqrt{\rho_{+}}}=\sigma.
\eqno{(4.8)}
$$
From   (4.2)  and (4.7),   we can easily get
$$\hspace{-18mm}\lim\limits_{k_{1}, k_{2}\rightarrow 0}\sigma_{2}=\lim\limits_{k_{1}, k_{2}\rightarrow 0}u_{\ast}+\lim\limits_{k_{1}, k_{2}\rightarrow 0}
             \sqrt{\frac{k_{1}}{\rho^{2}_{\ast}}
             +\frac{k_{2}^{2}\rho_{+}^{2}(\frac{1}{\rho_{\ast}}+\frac{1}{\rho_{+}})}{2\mu}}
             =\sigma.
\eqno{(4.9)}
$$
Thus it can be seen from (4.8) and (4.9) that  when $k_{1}, k_{2}\rightarrow 0, $  the  two shocks $S_{1}$ and $S_{2}$  will coincide whose velocities are identical with that  of the delta shock wave of  the transport equations with the same Riemann  initial data  $(\rho_{\pm}, u_{\pm})$.

Using the Rankine-Hugoniot conditions (2.7) for    $S_{1}$ and $S_{2}$, we have
$$
\left\{
             \begin{array}{ll}
               \sigma_{1}(\rho_{\ast}-\rho_{-})=\rho_{\ast}u_{\ast}-\rho_{-}u_{-}, \\
\sigma_{2}(\rho_{+}-\rho_{\ast})=\rho_{+}u_{+}-\rho_{\ast}u_{\ast}.
  \end{array}
\right.
\eqno{(4.10)}  $$
Then from (4.8)  and (4.9) it follows  that
$$\lim\limits_{k_{1}, k_{2}\rightarrow 0}(\sigma_{1}-
\sigma_{2})\rho_{\ast}=\lim\limits_{k_{1}, k_{2}\rightarrow 0}(\rho_{+}u_{+}-\rho_{-}u_{-}+\sigma_{1}\rho_{-}-\sigma_{2}\rho_{+})=[\rho u]-\sigma [\rho].\eqno{(4.11)}
$$
This  implies that
$$\lim\limits_{k_{1}, k_{2}\rightarrow 0}\int^{\sigma_{2}t}_{\sigma_{1}t}\rho_{\ast}dx=(\sigma[\rho]-[\rho u])t=\sqrt{\rho_{-}\rho_{+}}(u_{-}-u_{+})t=w(t).
\eqno{(4.12)}
$$

The proof is completed.

\vskip 0.1in

\vskip 0.1in

\noindent{\small {\small\bf Remark 1.} From the above results, it can be seen  that the limit of the Riemann solution  of system  (1.1) and (1.2)  as $k_{1}, k_{2}\rightarrow 0$ in the case  $(\rho_{+}, u_{+})\in IV(\rho_{-}, u_{-})$  is just the delta shock solution of (3.1)-(3.2)  when   $u_{-}> u_{+}$.
$$$$
\baselineskip 15pt
\sec{\large\bf 4.2.\quad   $\delta$-shocks  and concentration
}

Now, we  give the following results which give a very nice depiction of the limit in the case $(\rho_{+}, u_{+})\in IV(\rho_{-}, u_{-})$.
\vskip 0.1in
\noindent{\small {\small\bf Theorem 2.}    Let $u_{-}> u_{+}$  and $(\rho_{+}, u_{+})\in IV(\rho_{-}, u_{-})$. For any fixed $k_{1}, k_{2}>0,$  assuming that $(\rho, u)$ is a solution containing two  shocks $S_{1}$ and $S_{2}$
 of (1.1)-(1.2) with Riemann initial data (2.1), constructed in Section 2, it is obtained that as $k_{1}, k_{2}\rightarrow 0,$  $(\rho, u)$ converges in the sense of distributions, and the limit functions $\rho$ and $\rho u $ are the sums of a step function and a $\delta$-measure with weights
$$ (\sigma[\rho]-[\rho u])t\,\, \,{\mathrm\,and }\,\,\,\,(\sigma[\rho u]-[\rho u^{2}])t,
$$
respectively, which form a delta shock wave solution of (3.1) with the same Riemann  initial data  $(\rho_{\pm}, u_{\pm})$.
\vskip 0.1in
\noindent{\small\small\bf Proof.} Let $\xi= x/t$. Then for any fixed $k_{1}>0  $   and $ k_{2}>0,$ the Riemann solution to   the isentropic  magnetogasdynamics   equations  for  Chaplygin gas (1.1)-(1.2) can be written as
$$\hspace{5mm}(\rho, u)(\xi)=\left\{
                  \begin{array}{ll}
                    (\rho_{-}, u_{-}), & \hbox{$\xi<\sigma_{1}  $,} \\
                    (\rho_{\ast}, u_{\ast}), & \hbox{$\sigma_{1}<\xi<\sigma_{2}$,} \\
                    (\rho_{+}, u_{+}), & \hbox{$\xi>\sigma_{2} $,}
                  \end{array}
                \right.
\eqno{(4.13)}$$
which satisfies the following weak formulations:
$$\int_{-\infty}^{+\infty}(\xi-u(\xi))\rho(\xi)\psi'(\xi)d\xi
+\int_{-\infty}^{+\infty}\rho(\xi)\psi(\xi)d\xi=0\eqno{(4.14)}  $$
and
$$\int_{-\infty}^{+\infty}(\xi-u(\xi))\rho(\xi)u(\xi)\psi'(\xi)d\xi
-\int_{-\infty}^{+\infty}\Big(-\frac{k_{1}}{\rho(\xi)}
             +\frac{k_{2}^{2}(\rho(\xi))^{2}}{2\mu}\Big)
             \psi'(\xi)d\xi+\int_{-\infty}^{+\infty}\rho(\xi)u(\xi)\psi(\xi)d\xi=0\eqno{(4.15)}  $$
for any test function $\psi\in C_{0}^{\infty}(-\infty,  +\infty)$.

The first integral on the left-hand side of  (4.15) can be decomposed into
$$\bigg\{\int_{-\infty}^{\sigma_{1}}+\int_{\sigma_{1}}^{\sigma_{2}}+\int^{+\infty}_{\sigma_{2}}\bigg\}(\xi-u(\xi))
\rho(\xi)u(\xi )\psi'(\xi) d\xi. \eqno{(4.16)}  $$
The sum of the first and the last terms in (4.16) is
$$\int_{-\infty}^{\sigma_{1}}(\xi-u(\xi))
\rho(\xi)u(\xi )\psi'(\xi) d\xi +\int^{+\infty}_{\sigma_{2}}(\xi-u(\xi))
\rho(\xi)u(\xi )\psi'(\xi) d\xi$$$$ =\rho_{-}u_{-}\sigma_{1}\psi(\sigma_{1})-\rho_{-}u_{-}^{2}\psi(\sigma_{1})
-\rho_{-}u_{-}\int_{-\infty}^{\sigma_{1}}\psi(\xi) d\xi$$$$ -
\rho_{+}u_{+}\sigma_{2}\psi(\sigma_{2})+\rho_{+}u_{+}^{2}\psi(\sigma_{2})
-\rho_{+}u_{+}\int^{+\infty}_{\sigma_{2}}\psi(\xi) d\xi . \eqno{(4.17)}  $$
Taking the limit $k_{1}, k_{2}\rightarrow 0$  in (4.17) leads to
$$\hspace{-0mm}\lim\limits_{k_{1}, k_{2}\rightarrow 0}\bigg(\int_{-\infty}^{\sigma_{1}} +\int^{+\infty}_{\sigma_{2}}\bigg)(\xi-u(\xi))
\rho(\xi)u(\xi )\psi'(\xi) d\xi$$$$=([\rho u^{2}]-\sigma[\rho u])\psi(\sigma)-\int^{+\infty}_{-\infty}(\rho_{0}u_{0})(\xi-\sigma)\cdot\psi(\xi)d\xi,
\eqno{(4.18)}
$$
where $(\rho_{0}u_{0})(\xi)= \rho_{-}u_{-}+[\rho u]H(\xi)$ and $H $ is the Heaviside function.

For the second term in (4.16), integrating by parts  again,  we obtain
$$\hspace{-58mm}\int_{\sigma_{1}}^{\sigma_{2}}(\xi-u(\xi))
\rho(\xi)u(\xi )\psi'(\xi) d\xi=\int_{\sigma_{1}}^{\sigma_{2}}(\xi-u_{\ast})
\rho_{\ast}u_{\ast}\psi'(\xi) d\xi$$
$$\hspace{-26mm}=-\rho_{\ast}u_{\ast}^{2}(\psi(\sigma_{2})-\psi(\sigma_{1}))
+\rho_{\ast}u_{\ast}(\sigma_{2}\psi(\sigma_{2})-\sigma_{1}\psi(\sigma_{1}))-\rho_{\ast}u_{\ast}\int_{\sigma_{1}}^{\sigma_{2}}
\psi(\xi) d\xi  $$
$$=-u_{\ast}\rho_{\ast}(\sigma_{2}-\sigma_{1})\bigg(\frac{\psi(\sigma_{2})-\psi(\sigma_{1})}{\sigma_{2}-\sigma_{1}}u_{\ast}
-\frac{\sigma_{2}\psi(\sigma_{2})-\sigma_{1}\psi(\sigma_{1})}{\sigma_{2}-\sigma_{1}}+\frac{1}{\sigma_{2}-\sigma_{1}}
\int_{\sigma_{1}}^{\sigma_{2}}
\psi(\xi) d\xi\bigg). \eqno{(4.19)}  $$
Taking the limit $k_{1}, k_{2}\rightarrow 0$  in (4.19), noting (4.11) and the fact that both $\psi\in C_{0}^{\infty}(-\infty, +\infty)$ and $\lim\limits_{k_{1}, k_{2}\rightarrow 0}u_{\ast}=\lim\limits_{k_{1}, k_{2}\rightarrow 0}\sigma_{1}=\lim\limits_{k_{1}, k_{2}\rightarrow 0}\sigma_{2}=\sigma,$  we deduce that
$$\lim\limits_{k_{1}, k_{2}\rightarrow 0}\int_{\sigma_{1}}^{\sigma_{2}}(\xi-u(\xi))
\rho(\xi)u(\xi )\psi'(\xi) d\xi=\sigma([\rho u]-\sigma[\rho])(\sigma\psi'(\sigma)-\sigma\psi'(\sigma)-\psi(\sigma)+\psi(\sigma))
=0.\eqno{(4.20)}  $$

Similarly,  the first integral on the left-hand side of  (4.15)    can be decomposed into  three parts as
$$
-\bigg\{\int_{-\infty}^{\sigma_{1}}
             +\int_{\sigma_{1}}^{\sigma_{2}}+\int^{+\infty}_{\sigma_{2}}\bigg\}
\Big(-\frac{k_{1}}{\rho(\xi)}
             +\frac{k_{2}^{2}(\rho(\xi))^{2}}{2\mu}\Big)
             \psi'(\xi)d\xi,\eqno{(4.21)}  $$
which equals to
$$
\int_{-\infty}^{\sigma_{1}}
             \Big(\frac{k_{1}}{\rho_{-}}
             -\frac{k_{2}^{2}\rho_{-}^{2}}{2\mu}\Big)
             \psi'(\xi)d\xi
             +\int_{\sigma_{1}}^{\sigma_{2}}
\Big(\frac{k_{1}}{\rho_{\ast}}
             -\frac{k_{2}^{2}\rho_{\ast}^{2}}{2\mu}\Big)
             \psi'(\xi)d\xi
             +\int^{+\infty}_{\sigma_{2}}
\Big(\frac{k_{1}}{\rho_{+}}
             -\frac{k_{2}^{2}\rho_{+}^{2}}{2\mu}\Big)
             \psi'(\xi)d\xi$$
             $$=\Big(\frac{k_{1}}{\rho_{-}}
             -\frac{k_{2}^{2}\rho_{-}^{2}}{2\mu}\Big)\psi(\sigma_{1})
             +\frac{k_{1}}{\rho_{\ast}}(\psi(\sigma_{2})-\psi(\sigma_{1}))
-\frac{k_{2}^{2}\rho_{\ast}^{2}}{2\mu}(\psi(\sigma_{2})-\psi(\sigma_{1}))
-\Big(\frac{k_{1}}{\rho_{+}}
             -\frac{k_{2}^{2}\rho_{+}^{2}}{2\mu}\Big)\psi(\sigma_{2}).\eqno{(4.22)}  $$
Taking the limit $k_{1}, k_{2}\rightarrow 0$  in (4.22), by Lemmas 4.1-4.2,  we have
$$\hspace{-0mm}\lim\limits_{k_{1}, k_{2}\rightarrow 0}\int_{-\infty}^{+\infty}\Big(\frac{k_{1}}{\rho(\xi)}
             -\frac{k_{2}^{2}(\rho(\xi))^{2}}{2\mu}\Big)
             \psi'(\xi)d\xi=0.\eqno{(4.23)}  $$
Summarizing (4.18), (4.20)  and (4.23) leads to
$$\hspace{-0mm}\lim\limits_{k_{1}, k_{2}\rightarrow 0}\int_{-\infty}^{+\infty}((\rho u)(\xi)-(\rho_{0}u_{0})(\xi-\sigma))\psi(\xi)d\xi=(\sigma[\rho u]-[\rho u^{2}])\psi(\sigma),\eqno{(4.24)}  $$
which is true for any  $\psi\in C_{0}^{\infty}(-\infty,  +\infty)$.

As done previously, we can obtain the limit for the first  integral on the left-hand side of  (4.14) as
$$\hspace{-15mm}\lim\limits_{k_{1}, k_{2}\rightarrow 0}\int_{-\infty}^{+\infty}(\xi-u(\xi))\rho(\xi)\psi'(\xi)d\xi
=([\rho u]-\sigma[\rho])\psi(\sigma)
-\int_{-\infty}^{\sigma}\rho_{-}\psi(\xi)d\xi
-\int^{+\infty}_{\sigma}\rho_{+}\psi(\xi)d\xi$$$$\hspace{21mm}=([\rho u]-\sigma[\rho])\psi(\sigma)-\int^{+\infty}_{-\infty}\rho_{0}(\xi-\sigma)\psi(\xi)d\xi,\eqno{(4.25)}  $$
where $\rho_{0}(\xi)= \rho_{-}+[\rho ]H(\xi).$  Then returning to the formulation (4.14), we have
$$\hspace{-0mm}\lim\limits_{k_{1}, k_{2}\rightarrow 0}
\int_{-\infty}^{+\infty}(\rho(\xi)-\rho_{0}(\xi-\sigma))\psi(\xi)d\xi=(\sigma[\rho]-[\rho u])\psi(\sigma),\eqno{(4.26)}  $$
which is true for any  $\psi\in C_{0}^{\infty}(-\infty,  +\infty)$.

Finally, we  study the limits of $\rho$
 and $\rho u$ as $k_{1}, k_{2}\rightarrow 0$,  by tracing the time-dependence of weights of the $\delta$-measure.
Let $\phi(x, t)\in C_{0}^{\infty}((-\infty, +\infty)\times[0, +\infty))$, then we have
$$\hspace{-0mm}\lim\limits_{k_{1}, k_{2}\rightarrow 0}\int_{0}^{+\infty}\int_{-\infty}^{+\infty}\rho(x/t)\phi(x,t)dxdt=\lim\limits_{k_{1}, k_{2}\rightarrow 0}\int_{0}^{+\infty}t\bigg(\int_{-\infty}^{+\infty}\rho(\xi)\phi(\xi t,t)d\xi\bigg)dt.
\eqno{(4.27)}
$$
Regarding  $t$ as a parameter and applying (4.26),    one can easily see that
$$\hspace{-0mm}\lim\limits_{k_{1}, k_{2}\rightarrow 0}\int_{-\infty}^{+\infty}\rho(\xi)\phi(\xi t,t)d\xi=
\int_{-\infty}^{+\infty}\rho_{0}(\xi-\sigma)\phi(\xi t,t)d\xi+(\sigma[\rho]-[\rho u])\phi(\sigma t, t)
$$$$=\frac{1}{t}
\int_{-\infty}^{+\infty}\rho_{0}\Big(\frac{x}{t}-\sigma\Big)\phi(x,t)dx+(\sigma[\rho]-[\rho u])\phi(\sigma t, t),\eqno{(4.28)}  $$
Substituting (4.28) into (4.27)  and noting $\rho_{0}\big(\frac{x}{t}-\sigma\big)=\rho_{0}(x-\sigma t)$, we have
$$\hspace{-0mm}\lim\limits_{k_{1}, k_{2}\rightarrow 0}\int_{0}^{+\infty}\int_{-\infty}^{+\infty}\rho(x/t)\phi(x,t)dxdt= \int_{0}^{+\infty}\int_{-\infty}^{+\infty}\rho_{0}(x-\sigma t)\phi(x,t)dxdt$$$$+\int_{0}^{+\infty}t(\sigma[\rho]-[\rho u])\phi(\sigma t, t)dt.\eqno{(4.29)}  $$
By definition  (3.9), the last term on the right-hand side of (4.29) equals to $\langle w_{1}(t)\delta_{S},\phi(\cdot, \cdot)\rangle$, where
$$w_{1}(t)=  (\sigma[\rho]-[\rho u])t.$$
With the same reason as before, we arrive at
$$\hspace{-0mm}\lim\limits_{k_{1}, k_{2}\rightarrow 0}\int_{0}^{+\infty}\int_{-\infty}^{+\infty}\rho(x/t)u(x/t)\phi(x,t)dxdt= \int_{0}^{+\infty}\int_{-\infty}^{+\infty}(\rho_{0}u_{0})(x-\sigma t)\phi(x,t)dxdt$$$$+\int_{0}^{+\infty}t(\sigma[\rho u]-[\rho u^{2}])\phi(\sigma t, t)dt.\eqno{(4.30)}  $$
The last term on the right-hand side of (4.30) equals to $\langle w_{2}(t)\delta_{S},\phi(\cdot, \cdot)\rangle$, where
$$w_{2}(t)=  (\sigma[\rho u]-[\rho u^{2}])t.$$
The proof is completed.
\vskip 0.1in

\baselineskip 15pt
 \sec{\Large\bf 5.\quad  Formation of vacuum states}
In this section,  we study the formation of  vacuum states in the  Riemann solutions of system  (1.1) and (1.2) in the case  $(\rho_{+}, u_{+})\in I(\rho_{-}, u_{-})$  with  $u_{-}<u_{+}$ and $\rho_{\pm}>0$ as  both the pressure and the magnetic field vanish.
In this case,  we know that the Riemann solution   consists of a backward rarefaction wave $R_{1}$, a
forward rarefaction wave $R_{2}$  and an intermediate state $(\rho_{\ast}, u_{\ast})$  besides two constant states $(\rho_{\pm}, u_{\pm})$, which are as follows
$$R_{1}:
\left\{
  \begin{array}{ll}
    \xi=\lambda_{1}=u-\sqrt{\frac{k_{1}}{\rho^{2}}+\frac{k_{2}^{2}\rho}{\mu}}, \\
    u=u_{-}-\int^{\rho}_{\rho_{-}}\frac{\sqrt{\frac{k_{1}}{s^{2}}+\frac{k_{2}^{2}s}{\mu}}}{s}ds,\,\,\,\,
    \rho_{\ast}\leq\rho\leq\rho_{-},
  \end{array}
\right.
\eqno{(5.1)}  $$
and$$R_{2}:
\left\{
  \begin{array}{ll}
    \xi=\lambda_{2}=u+\sqrt{\frac{k_{1}}{\rho^{2}}+\frac{k_{2}^{2}\rho}{\mu}}, \\
    u=u_{+}+\int^{\rho}_{\rho_{+}}\frac{\sqrt{\frac{k_{1}}{s^{2}}+\frac{k_{2}^{2}s}{\mu}}}{s}ds,\,\,\,\,
    \rho_{\ast}\leq\rho\leq\rho_{+}.
  \end{array}
\right.
\eqno{(5.2)}  $$
From (5.1) and (5.2), we can  derive
 $$ u_{+}-u_{-}=\int^{\rho_{-}}_{\rho_{\ast}}\frac{\sqrt{\frac{k_{1}}{s^{2}}+\frac{k_{2}^{2}s}{\mu}}}{s}ds
 +\int_{\rho_{\ast}}^{\rho_{+}}\frac{\sqrt{\frac{k_{1}}{s^{2}}+\frac{k_{2}^{2}s}{\mu}}}{s}ds,\,\,\,\,
    \rho_{\ast}\leq\rho_{\pm}.  \eqno{(5.3)}
$$
For any given $\rho_{\pm}>0$, if $\lim\limits_{k_{1}, k_{2}\rightarrow 0}\rho_{\ast}=K\in(0, \min(\rho_{-}, \rho_{+})],$  then
by $$ \int_{\rho_{\ast}}^{\rho}\frac{\sqrt{A+\frac{B}{s^{2}}}}{s}ds=
 -\sqrt{A+\frac{B}{\rho^{2}}}
 +\sqrt{A}\ln\bigg(\sqrt{A+\frac{B}{\rho^{2}}}
 +\sqrt{A}\bigg)+\sqrt{A}\ln\rho
$$$$+\sqrt{A+\frac{B}{\rho_{\ast}^{2}}}
 -\sqrt{A}\ln\bigg(\sqrt{A+\frac{B}{\rho_{\ast}^{2}}}
 +\sqrt{A}\bigg)-\sqrt{A}\ln\rho_{\ast},\,\,\,\,\,A>0,\eqno{(5.4)}
$$it follows that
  $$ 0\leq\int^{\rho_{-}}_{\rho_{\ast}}\frac{\sqrt{\frac{k_{1}}{s^{2}}+\frac{k_{2}^{2}s}{\mu}}}{s}ds
 \leq\int_{\rho_{\ast}}^{\rho_{-}}\frac{\sqrt{\frac{k_{1}}{s^{2}}+\frac{k_{2}^{2}\rho_{-}}{\mu}}}{s}ds$$
 $$= \sqrt{\frac{k_{1}}{\rho_{\ast}^{2}}+\frac{k_{2}^{2}\rho_{-}}{\mu}}
 -\sqrt{\frac{k_{2}^{2}\rho_{-}}{\mu}}\ln\Bigg(\sqrt{\frac{k_{1}}{\rho_{\ast}^{2}}+\frac{k_{2}^{2}\rho_{-}}{\mu}}
 +\sqrt{\frac{k_{2}^{2}\rho_{-}}{\mu}}\Bigg)-\sqrt{\frac{k_{2}^{2}\rho_{-}}{\mu}}\ln\rho_{\ast}
    $$
$$-\sqrt{\frac{k_{1}}{\rho_{-}^{2}}+\frac{k_{2}^{2}\rho_{-}}{\mu}}
 +\sqrt{\frac{k_{2}^{2}\rho_{-}}{\mu}}\ln\Bigg(\sqrt{\frac{k_{1}}{\rho_{-}^{2}}+\frac{k_{2}^{2}\rho_{-}}{\mu}}
 +\sqrt{\frac{k_{2}^{2}\rho_{-}}{\mu}}\Bigg)+\sqrt{\frac{k_{2}^{2}\rho_{-}}{\mu}}\ln\rho_{-}\rightarrow 0,\,\,\,\,
   {\mathrm as\,\,  \,}  k_{1}, k_{2}\rightarrow 0.  \eqno{(5.5)}
$$
Therefore, by the squeeze theorem  in multivariable calculus,  we arrive at
$$\lim\limits_{k_{1}, k_{2}\rightarrow 0}\int^{\rho_{-}}_{\rho_{\ast}}\frac{\sqrt{\frac{k_{1}}{s^{2}}+\frac{k_{2}^{2}s}{\mu}}}{s}ds=0
.
\eqno{(5.6)}
$$
Similarly, we can obtain  that $$ \lim\limits_{k_{1}, k_{2}\rightarrow 0}\int_{\rho_{\ast}}^{\rho_{+}}\frac{\sqrt{\frac{k_{1}}{s^{2}}+\frac{k_{2}^{2}s}{\mu}}}{s}ds=0.  \eqno{(5.7)}
$$
Combining (5.3), (5.6) and (5.7), we have $u_{-}-u_{+}=0,$
which contradicts with $u_{-}<u_{+}.$
Therefore, $\lim\limits_{k_{1}, k_{2}\rightarrow 0}\rho_{\ast}=0$, which implies that a vacuum  occurs.
  From (5.1), one can  see that

 $$u_{-}-\sqrt{\frac{k_{1}}{\rho^{2}}+\frac{k_{2}^{2}\rho}{\mu}}\leq\lambda_{1}=u_{-}-\sqrt{\frac{k_{1}}{\rho^{2}}
 +\frac{k_{2}^{2}\rho}{\mu}}+\int^{\rho_{-}}_{\rho}\frac{\sqrt{\frac{k_{1}}{s^{2}}+\frac{k_{2}^{2}s}{\mu}}}{s}ds
 $$$$\leq u_{-}-\sqrt{\frac{k_{1}}{\rho^{2}}
 +\frac{k_{2}^{2}\rho}{\mu}}+\int^{\rho_{-}}_{\rho}\frac{\sqrt{\frac{k_{1}}{s^{2}}+\frac{k_{2}^{2}\rho_{-}}{\mu}}}{s}ds, \,\,\,\,\,
    \rho_{\ast}\leq\rho\leq\rho_{-} \eqno{(5.8)}$$
 It can be derived from (5.4) that
 $$ u_{-}-\sqrt{\frac{k_{1}}{\rho^{2}}+\frac{k_{2}^{2}\rho}{\mu}}+\int_{\rho}^{\rho_{-}}\frac{\sqrt{\frac{k_{1}}{s^{2}}+\frac{k_{2}^{2}\rho_{-}}{\mu}}}{s}ds$$
 $$=u_{-}-\sqrt{\frac{k_{1}}{\rho^{2}}+\frac{k_{2}^{2}\rho}{\mu}}+ \sqrt{\frac{k_{1}}{\rho^{2}}+\frac{k_{2}^{2}\rho_{-}}{\mu}}
 -\sqrt{\frac{k_{2}^{2}\rho_{-}}{\mu}}\ln\Bigg(\sqrt{\frac{k_{1}}{\rho^{2}}+\frac{k_{2}^{2}\rho_{-}}{\mu}}
 +\sqrt{\frac{k_{2}^{2}\rho_{-}}{\mu}}\Bigg)-\sqrt{\frac{k_{2}^{2}\rho_{-}}{\mu}}\ln\rho
    $$
$$-\sqrt{\frac{k_{1}}{\rho_{-}^{2}}+\frac{k_{2}^{2}\rho_{-}}{\mu}}
 +\sqrt{\frac{k_{2}^{2}\rho_{-}}{\mu}}\ln\Bigg(\sqrt{\frac{k_{1}}{\rho_{-}^{2}}+\frac{k_{2}^{2}\rho_{-}}{\mu}}
 +\sqrt{\frac{k_{2}^{2}\rho_{-}}{\mu}}\Bigg)+\sqrt{\frac{k_{2}^{2}\rho_{-}}{\mu}}\ln\rho_{-}.\,\,\,\,
     \eqno{(5.9)}
$$
  The uniform boundedness of $\rho(\xi)$ with respect to $k_{1}, k_{2}$  in this case leads to
$$\lim\limits_{k_{1}, k_{2}\rightarrow 0}\bigg(u_{-}-\sqrt{\frac{k_{1}}{\rho^{2}}+\frac{k_{2}^{2}\rho}{\mu}}
+\int_{\rho}^{\rho_{-}}\frac{\sqrt{\frac{k_{1}}{s^{2}}+\frac{k_{2}^{2}\rho_{-}}{\mu}}}{s}ds\bigg)
=\lim\limits_{k_{1}, k_{2}\rightarrow 0}\bigg(u_{-}-\sqrt{\frac{k_{1}}{\rho^{2}}+\frac{k_{2}^{2}\rho}{\mu}}
\bigg)=u_{-}.\eqno{(5.10)}  $$
Then, by the squeeze theorem  in multivariable calculus,  we have
$\lim\limits_{k_{1}, k_{2}\rightarrow 0}\lambda_{1}=u_{-}.
$
Similarly, we can obtain  that $$\lim\limits_{k_{1}, k_{2}\rightarrow 0}\lambda_{2}=u_{+}
\, \,{\mathrm\,and }\,\,\,
\lim\limits_{k_{1}, k_{2}\rightarrow 0}u(\xi)=\xi,\,\,\,\, \,{\mathrm\,for }\,\,\,  \xi\in(u_{-}, u_{+}).
\eqno{(5.11)}
$$
Then from above we have proved  the following results.
\vskip 0.1in
\noindent{\small {\small\bf Theorem 3.}  In the case  $(\rho_{+}, u_{+})\in I(\rho_{-}, u_{-})$  with  $u_{-}<u_{+},$ as $k_{1}, k_{2}\rightarrow 0$,  the vacuum  state occurs and two rarefaction waves $R_{1}$ and $R_{2}$ become two contact discontinuities
$u=u_{-}$ and $u=u_{+}$, respectively, connecting the constant states $(\rho_{\pm}, u_{\pm})$ with the vacuum $(\rho =0).$

\vskip 0.1in
\noindent{\small {\small\bf Theorem 4.} In the case  $(\rho_{+}, u_{+})\in I(\rho_{-}, u_{-})$  with  $u_{-}<u_{+},$ as $k_{1}, k_{2}\rightarrow 0$, the limit of the Riemann solution  of  (1.1) and  (1.2)  with initial data (2.1)  is just    the Riemann solution    of  the transport equations (3.1) for
zero pressure flow with the same  initial data, which contains  two contact discontinuities $\xi=x/t=u_{\pm}$  and a vacuum  state besides two constant states.

\end{document}